\documentclass[3p,times,procedia]{elsarticle}
\usepackage{ecrc}
\usepackage{amsmath}
\usepackage{algorithm}
\usepackage{algorithmicx}
\usepackage{subfig}
\usepackage{algpseudocode}
\usepackage[utf8]{inputenc}

\volume{00}

\firstpage{1}

\journalname{Transportation Research Procedia}

\runauth{M.A Rahman}

\jid{trpro}

\usepackage{amssymb}

\biboptions{authoryear}

\usepackage[figuresright]{rotating}

\usepackage[bookmarks=false]{hyperref}
    \hypersetup{colorlinks,
      linkcolor=blue,
      citecolor=blue,
      urlcolor=blue}

\begin{document}
\begin{frontmatter}



\dochead{The 1st International Conference on Smart Mobility and Logistics Ecosystems (SMiLE)
September 17-19, 2024, KFUPM, Saudi Arabia}%

\title{Logistics Hub Location Optimization: A K-Means and P-Median Model Hybrid Approach Using Road Network Distances}


\author[a]{M.A Rahman} 
\author[b,c]{M.A Basheer\corref{cor1}}
\author[a]{Z. Khalid}
\author[a]{M. Tahir}
\author[a]{M. Uppal}

\address[a]{Department of Electrical Engineering, LUMS, Lahore 54792, Pakistan}
\address[b]{Department of Architecture and City Design, King Fahd University of Petroleum and Minerals, Dhahran 5067, Saudi Arabia}
\address[c]{Interdisciplinary Center of Smart Mobility and Logistics, King Fahd University of Petroleum and Minerals, Dhahran 5067, Saudi Arabia}

\begin{abstract}
Logistic hubs play a pivotal role in the last-mile delivery distance; even a slight increment in distance negatively impacts the business of the e-commerce industry while also increasing its carbon footprint. The growth of this industry, particularly after Covid-19, has further intensified the need for optimized allocation of resources in an urban environment. In this study, we use a hybrid approach to optimize the placement of logistic hubs. The approach sequentially employs different techniques. Initially, delivery points are clustered using K-Means in relation to their spatial locations. The clustering method utilizes road network distances as opposed to Euclidean distances. Non-road network-based approaches have been avoided since they lead to erroneous and misleading results. Finally, hubs are located using the P-Median method. The P-Median method also incorporates the number of deliveries and population as weights. Real-world delivery data from Muller and Phipps (M\&P) is used to demonstrate the effectiveness of the approach. Serving deliveries from the optimal hub locations results in the saving of 815 (10\%) meters per delivery.
\end{abstract}

\begin{keyword}
Urban logistics; Optimization; Last mile delivery; Machine Learning; K-Means; P-Median




\end{keyword}
\cortext[cor1]{Corresponding author. Tel.: +966-546-325-954 }
\end{frontmatter}

\email{m.basheer@kfupm.edu.sa}



\section{Introduction}
\label{Introduction}
The need for urban logistics has been expanding quickly in terms of both scope and intensity throughout the current economic globalization and e-commerce boom (\citeauthor{bib1}, \citeyear{bib1}). The Covid-19 pandemic has further accelerated this trend by changing people's traditional buying habits, as observed by the improved financial performance of logistics organizations in 2020 (\citeauthor{bib2}, \citeyear{bib2}). For instance, the e-commerce sector generated \$3.351 trillion in global sales in 2019 and is predicted to reach \$6.169 trillion by 2023 (\citeauthor{bib3}, \citeyear{bib3}). Because of this upward trend, small, medium-sized, and micro businesses have been able to bear the brunt of Covid-19. Similar growth has been observed in Pakistan, where the value of the e-commerce sector was \$548.89 million in 2021 with a growth of 35\% in the first quarter (\citeauthor{bib4}, \citeyear{bib4}).

This growth has consequently led to more delivery parcels being packed, shipped, and transported by logistics providers. All parcels invariably pass through a Logistics Hub (LH) on their way to the customer. A LH may be defined as a distribution center that connects other transportation systems including local, national, and international inside which parcels may be sorted, consolidated, stored, and finally shipped to their respective locations (\citeauthor{bib5}, \citeyear{bib5}). The logistics system as a whole is affected by the decision of where to locate facilities across the network. The unoptimized location of a LH leads to delivery vehicles traveling longer distances hence leading to negative spillovers such as emission of harmful gasses into the environment. As a result, a LH has a significant effect on carbon emissions in the delivery process \cite{bib6} since logistics transportation considerably contributes to environmental pollution and greenhouse gas effects (\citeauthor{bib7}, \citeyear{bib7}). For example, road freight transport emitted over 1.8 billion tonnes of CO$_2$ in 2019 (\citeauthor{bib8}, \citeyear{bib8}). Spillovers also include higher delivery costs. Last-mile delivery accounts for nearly $53\%$ of total shipping costs (\citeauthor{bib9}, \citeyear{bib9}). Fuel costs account for approximately $10\%$ to $25\%$ percent of the overall last-mile delivery expenses (\citeauthor{bib10}, \citeyear{bib10}). Traveling longer distances due to unoptimized hub locations leads to more fuel consumption and therefore higher total costs for last-mile deliveries. Moreover, our discussion with logistics operators also highlighted the location of LH as one of the most significant problems to solve. 

The hub location problem is a significant issue with several applications in fields like telecommunications and transportation (\citeauthor{bib11}, \citeyear{bib11}). One such domain is the logistics industry. Most studies that try to solve this Multi-Facility Location Problem (MFLP) use non-road network distances. Euclidean or other sorts of distance calculations are used (\citeauthor{bib12}, \citeyear{bib12}). Numerous optimization techniques have been used to tackle the problem of hub location optimization. \cite{bib5} use multiple criteria techniques for selecting suitable locations for LH. The novel Mixed Integer Quadratic Constrained Programming (MIQCP) method is used to cluster centroids representing the most suitable sites. A weighted K-Means algorithm has also been employed in studies to determine the optimal location of container storage and distribution centers (\citeauthor{bib13}, \citeyear{bib13}; \citeauthor{bib14}, \citeyear{bib14}). Multiple characteristics of potential locations as weights are used to determine optimal points. \cite{bib15} propose a two-part method. First making use of K-Means clustering after which in the second part, for each cluster a Mixed Integer Linear Programming (MILP) technique is implemented to help determine the maximum profit-yielding facility. \cite{bib16} utilize a Fuzzy clustering-based hybrid method to locate multiple facilities. Firstly, customers are grouped based on geographical location using cluster analysis. Then, cluster centers are used to locate facilities. Finally, every cluster is solved independently, treating each as a single facility location problem. The center of gravity method is used to improve facility location. Similarly, \cite{bib17} apply fuzzy analytic hierarchy process, spatial statistics, and analysis approaches to determine logistics center locations.
\cite{bib18} use a two-part methodology to optimize hub placement. In the first part, the Axiomatic Fuzzy Set (AFS) clustering method is employed to evaluate the logistics center location. In the second part, a Technique for Order Preference by Similarity to Ideal Solution (TOPSIS) is used to determine the final selection. \cite{bib19} use a Multi-Criteria Decision Making (MCDM) model built on a Geographic Information System (GIS) to assess suitable sites for freight communities. A four-step model is used to define weight through the Analytical Network Process (ANP) and to find the optimal location using the TOPSIS technique. \cite{bib20} focus on Intermodal Hub-and-Spoke (IH\&S) network design and optimization taking into account mixed uncertainties in both transportation cost and trip time. \cite{bib21} employ fuzzy SWARA and CoCoSo methods to locate logistics centers. However, the distance metric used for optimization in all these studies is based on non-road network distances, with Euclidean calculations being the most common. Euclidean distance is defined as the straight line distance between two points in a plane \cite{bib22}, underestimating road network distance. Non-road network-based approaches could lead to erroneous and misleading results since they misrepresent ground distances. Some researchers have tried to overcome this issue by multiplying Euclidean calculations with a constant factor of 1.25 to transform them into real-world distances (\citeauthor{bib16}, \citeyear{bib16}). Similarly, \cite{bib13} use a novel optimization method to determine the location of the logistics hub which is closer to the most distant of the roads. Roads are treated as straight lines, and the algorithm uses custom, non-road network distances. \cite{bib23} implement a variation of Euclidean distance for cluster analysis. In this approach, weighted spatial patterns are used to find optimal locations for LH. For distance metric, Euclidean is used from centroid to the point. The distance is normalized by the sum of distances from all centroids. Nevertheless, these approaches do not incorporate true network distance for LH optimization.

In the literature, road networks have been used widely for station location and other problems to determine the most optimal site. \cite{bib24} use road network distance for optimal dynamic traffic network planning. \cite{bib25} use road networks as input to determine the optimal set of locations in the supply chain of a firm that is operating under competition. \cite{bib26} use road networks with average speeds to determine the optimal location of hydrogen stations in Florida. However, a road network distance matrix has not been employed to solve the problem of location optimization for LH. The contributions of this study are as follows:
\begin{enumerate}
\item Incorporate ground truth road network for computation of accurate distances in the location optimization problem.
\item Propose a hybrid approach utilizing both K-Means and P-Median to determine optimal hub locations serving the deliveries while minimizing travel distance.
\item Demonstrate the effectiveness of the proposed approach using real-world delivery data.
\end{enumerate}
The rest of the paper is structured as follows: Methodology (section 2) explains the hybrid approach using road network distances, Application and Analysis (section 3) outlines our case study area and data sources used including baseline measures, Results (section 4) demonstrates effectiveness of the proposed approach and Conclusion (section 5) provides concluding remarks.
\section{Methodology}\label{Methodology}

\subsection{Problem Description and Formulation}
The proposed approach attempts to answer the following question: \emph{Given $M$ number of LH clusters defined as a set of spatial points and  $T$ total delivery points, where should hubs be placed such that average road network distance from hubs to the respective locations per delivery is minimized?} This can be formulated as the following optimization problem with an objective function minimizing the average distance traveled per delivery:
\begin{equation}\label{equation 1}
{\rm minimize}\quad \frac{\sum\limits_{j=0}^{M-1} \sum\limits_{i=1}^{T} a_{ij} D(x_i,c_j)}{T} ,
\end{equation}
where $x_i$ is unique delivery location and $a_{ij}$ = 1 if $x_i$ $\in$ $j$ cluster, else $a_{ij}$ = 0. The function $D$ returns road network distance from hub $c_j$ of cluster $j$ to point $x_i$. Distance between points is calculated using a road network with directions incorporated indicating whether a road is two-way or one-way. 
\subsection{Cluster Formation}
To optimize the above objective function, the proposed approach makes use of the K-Means algorithm which follows two steps. In the first step, the distance between all the points and initial centroids is calculated. The initial centroids can be any random spatial points and not necessarily LH. Based on these distances, points are associated with the nearest centroid. In the second step, based on the new association of points with the nearest centroid, the location of the centroids of each cluster is updated. These steps are repeated until the maximum number of iterations is reached or changes in the location of the centroids are minimal. The maximum number of iterations and change cutoff are input to the algorithm (\citeauthor{bib27}, \citeyear{bib27}). The distances in the first step are calculated using the QNEAT3 plugin in QGIS which makes use of road network. The plugin uses Dijkstra's algorithm to determine the shortest distance between origin and destination (\citeauthor{bib28}, \citeyear{bib28}). 

\subsection{Hub Location Optimization}
Traditionally, the centroid update step (second step) uses the mean position of points in the cluster to update centroids. However, this approach is not suitable since centroids/LH need to be updated based on road network distances. To account for this point, we use the well-defined P-Median model as part of the algorithm. The P-Median model is defined as follows:
\begin{align}
\quad & {\rm minimize} \quad \sum_{j=0}^{J} \sum_{i=0}^{I} h_i d_{ij} Y_{ij} \\ 
\textrm{such that} \quad & \sum_{j=0}^{J} Y_{ij} = 1 \quad \forall \, i = 0,1,2...I,  \\
    & Y_{ij} - X_j \leq 0 \quad \forall \, i, \forall \, j = 0,1,2...J, \\
    &  \sum_{j=0}^{J} X_j = p,
\end{align} 
where $J$ is the total possible hub locations , $I$ is the total number of unique delivery locations, $h_i$ is the weight of delivery location, $d_{ij}$ quantifies road network distance between $i$-th delivery location and $j$-th hub location, $Y_{ij}$ is 1 if $i$ is assigned to $j$ else 0, $X_j$ is 1 if LH is located at $j$ else 0 and $p$ is the number of LH to be determined.  The first constraint ensures that an $i$ is assigned to only one $j$. The second constraint makes sure that an $i$ is assigned to $j$ only if a LH exists at $j$. The third constraint is there to make certain that the total number of LH is $p$. Since this approach locates only 1 center when updating the centroid for each cluster, $p$ is set to 1. To determine $J$, the cluster is divided into a grid of resolution 1 km x 1 km. The center of the grid is attached to the nearest road network node. Using these as possible hub locations, the P-Median problem mentioned above is solved using the QNEAT3 plugin inputting road network distances with directions. For weight $h_i$, the number of deliveries performed at unique location $i$ is used in phase 1. In phase 2, the above methodology is repeated using a layer of population estimates to make the approach more robust. For this purpose, the WorldPop data set is used that provides population estimates at a 1km x 1km resolution. Each pixel in the dataset provides a value of population i.e., population per pixel (ppp) (\citeauthor{bib34}, \citeyear{bib34}). Using this data set, pixels are converted to 1km x 1km polygons, and their centroids are determined. These points are used to form clusters and determine LH. Weight of each point, $h_i$ is calculated as follows:
\begin{equation}\label{equation 1}
h_i = \alpha x + (1-\alpha) y,
\end{equation}
where $x$ is the number of normalized deliveries performed and $y$ is normalized ppp in grid $i$. $\alpha$ is a hyper-parameter quantifying the weight given to each term in $h_i$. 







\begin{figure}
\centering
\subfloat[Iteration 1]{
  \includegraphics[width=0.45\columnwidth, trim={0 2cm 0 2cm}, clip]{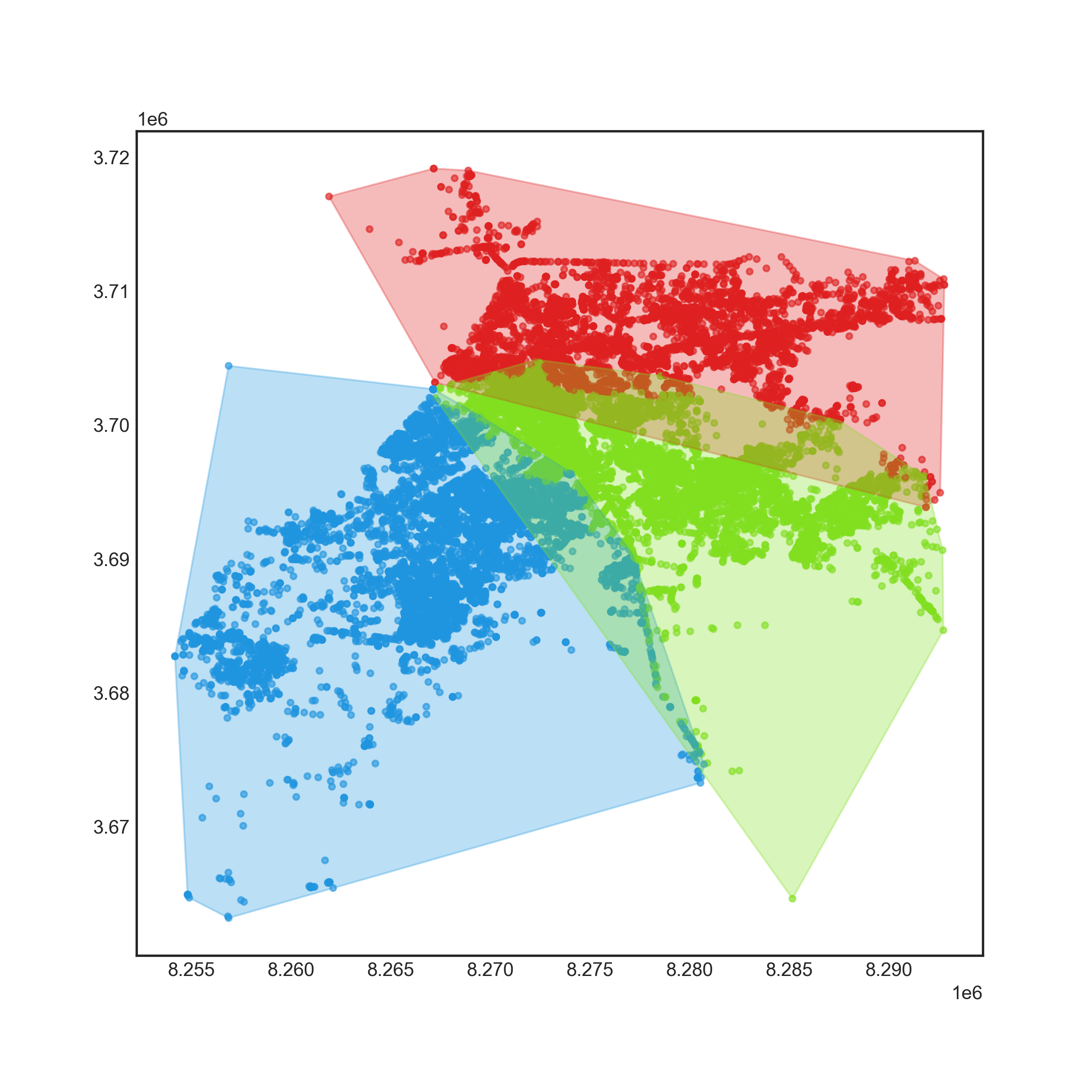}
  \label{Fig: 1(a)}
}
\subfloat[Iteration 2]{
  \includegraphics[width=0.45\columnwidth, trim={0 2cm 0 2cm}, clip]{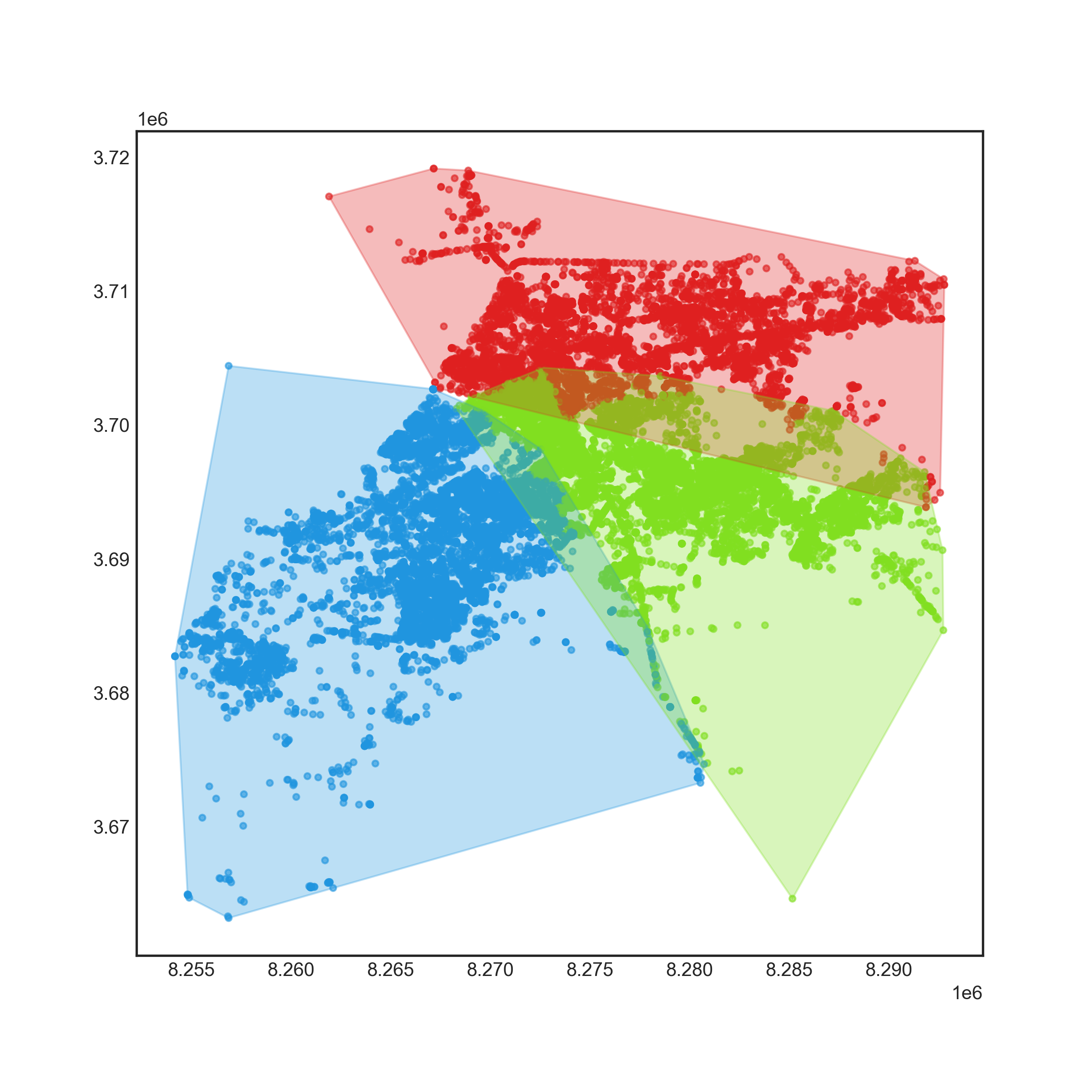}
  \label{Fig: 1(b)}
}
\vfill
\subfloat[Iteration 3]{
  \includegraphics[width=0.45\columnwidth, trim={0 2cm 0 2cm}, clip]{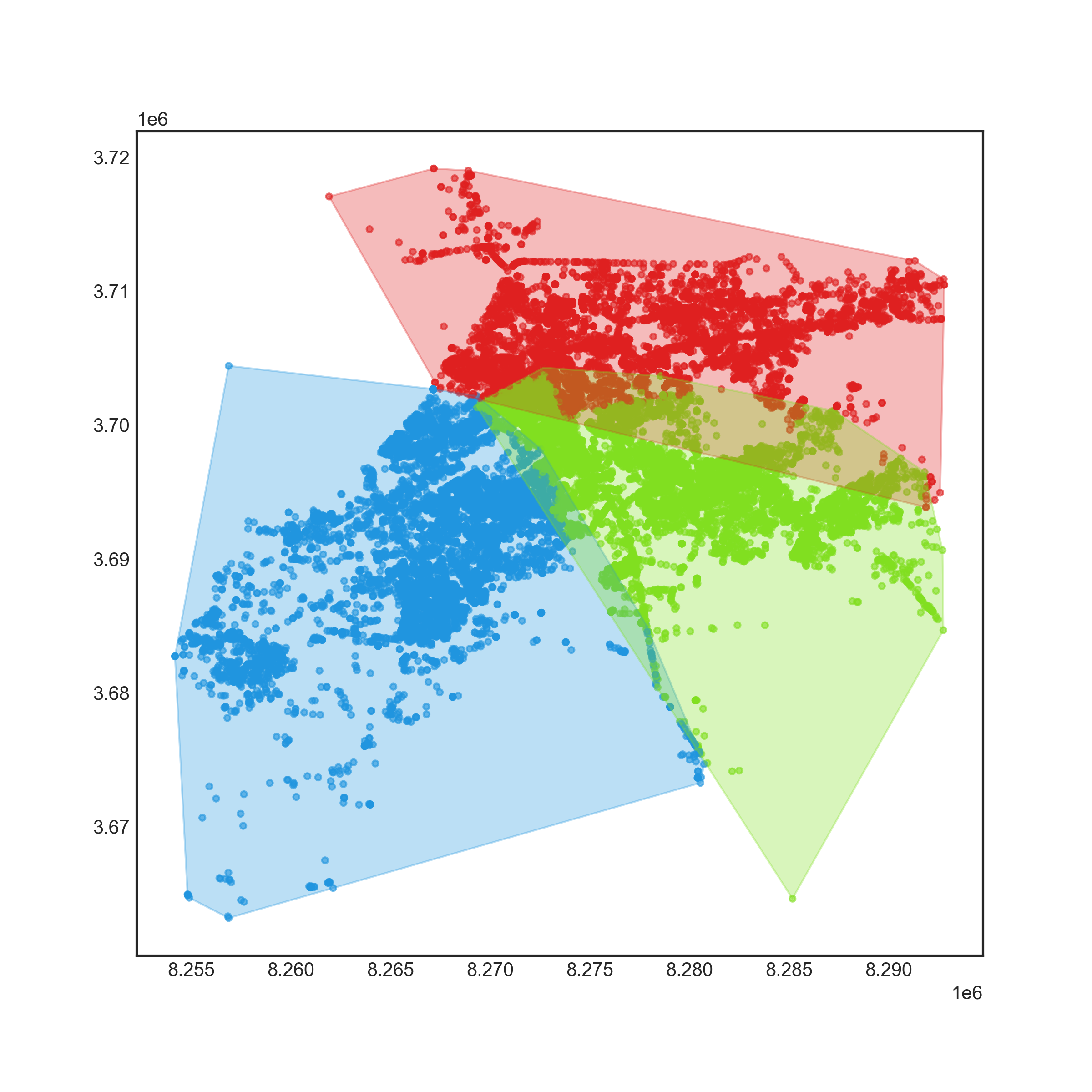}
  \label{Fig: 1(c)}
}
\subfloat[Iteration 4]{
  \includegraphics[width=0.45\columnwidth, trim={0 2cm 0 2cm}, clip]{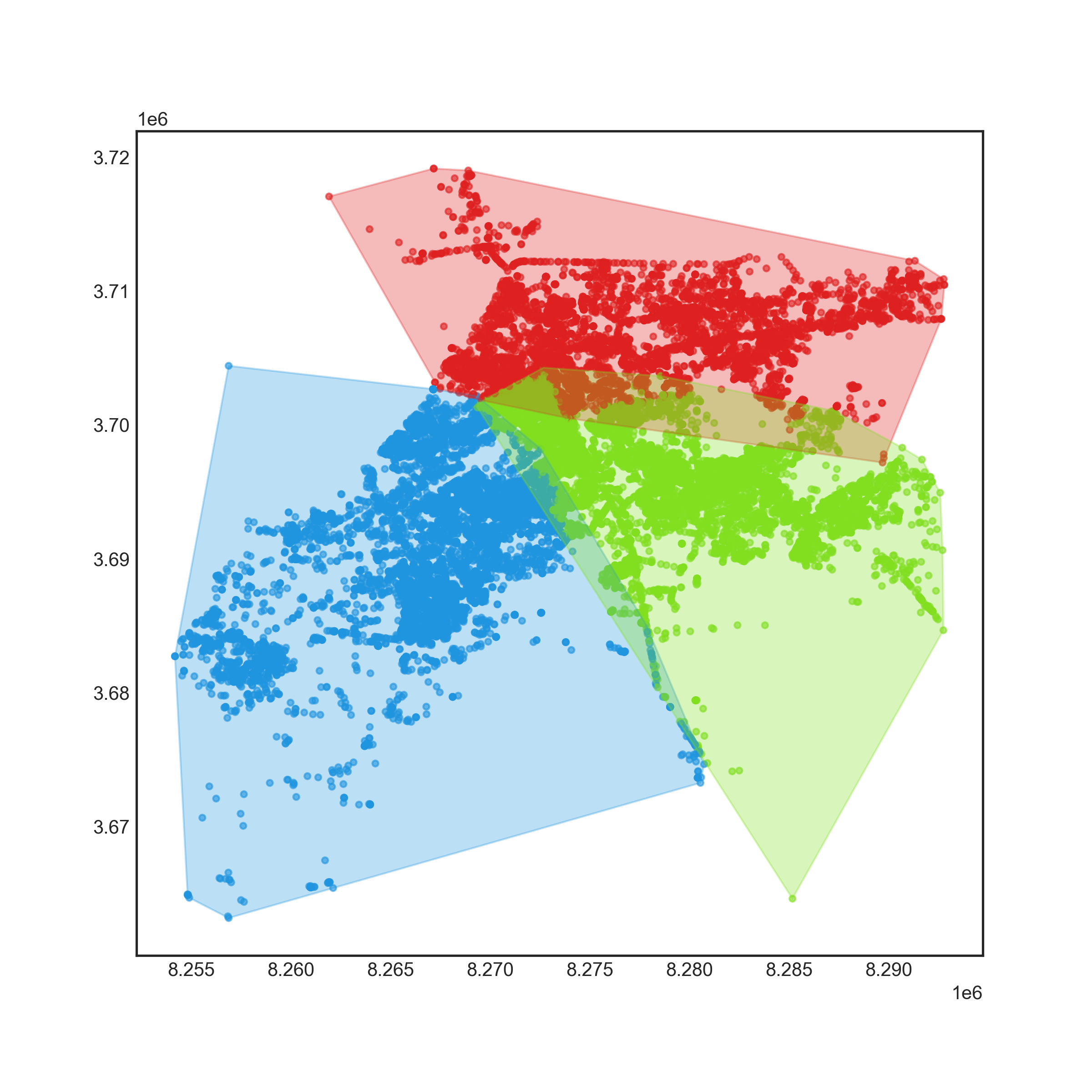}
  \label{Fig: 1(d)}
}
\vfill
\subfloat[Iteration 5]{
  \includegraphics[width=0.45\columnwidth, trim={0 2cm 0 2cm}, clip]{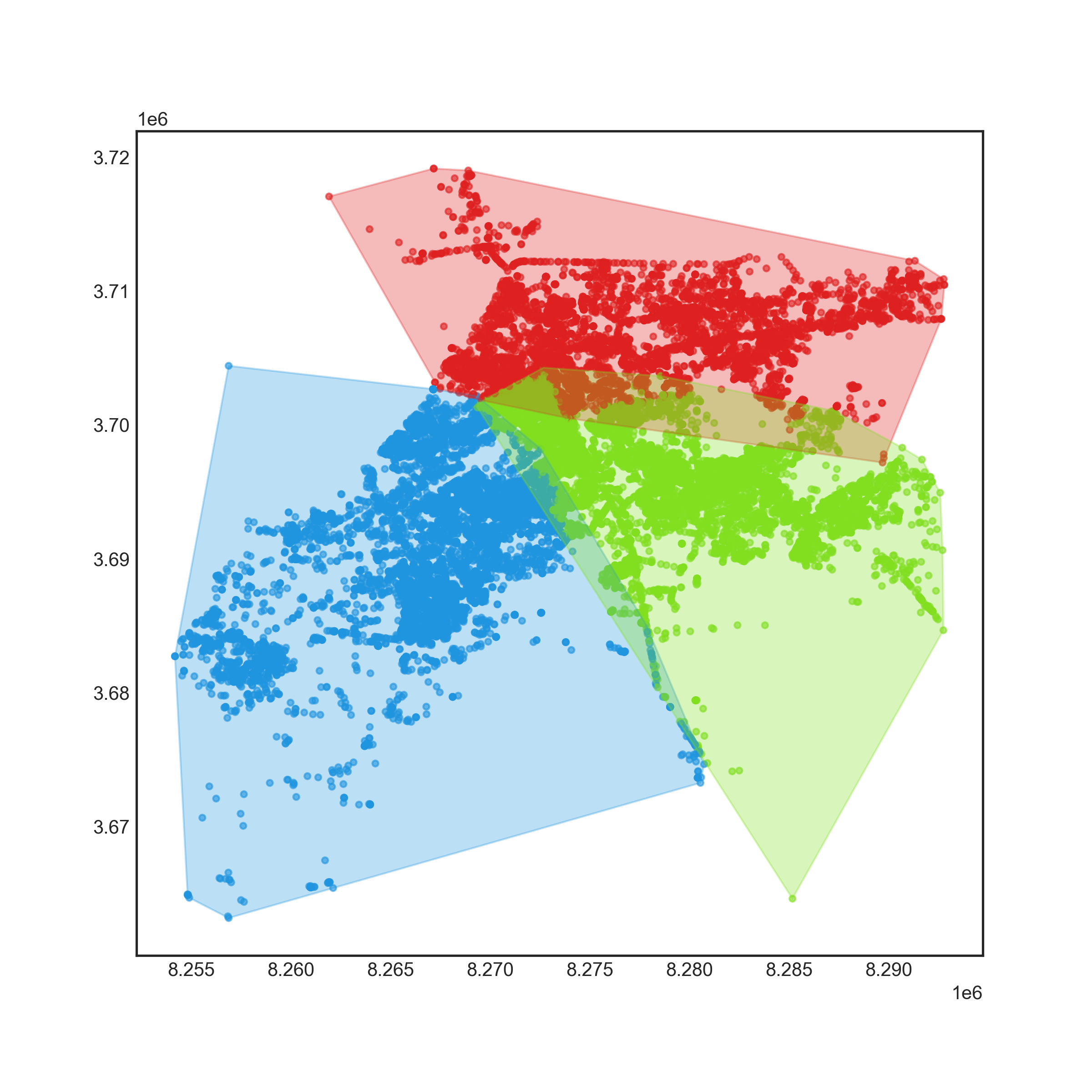}
  \label{Fig: 1(e)}
}
\subfloat[Average distance per delivery vs iteration number using M\&P Data for phase 1]{
  \includegraphics[width=0.45\columnwidth, trim={0 2cm 0 2cm}, clip]{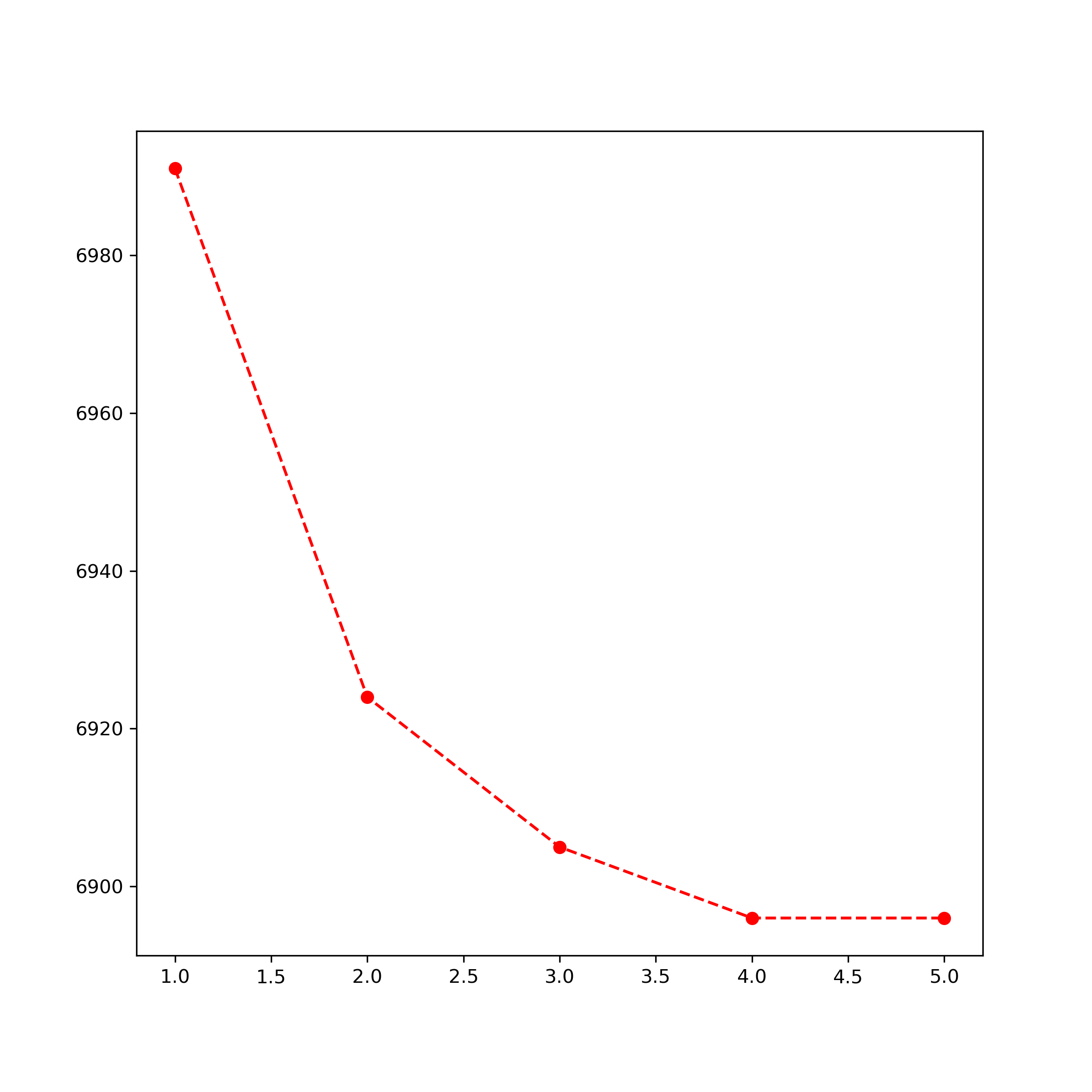}
  \label{Fig: 1(f)}
}
\caption{(a) - (e) Convex hull of each cluster formed from iteration 1 to 5. (f) Monotonically decreasing curve showing the change in average distance per delivery with each iteration.}
\label{Fig: 1}
\end{figure}


\section{Application and Analysis}\label{Case Study}
This study takes Lahore, Pakistan as a case to draw inferences about LH optimization. The city has a population of around 11.13 million, with approximately 1,744,755 households (\citeauthor{bib29}, \citeyear{bib29}). The GDP of the city is estimated to be over 1 trillion rupees \cite{bib4} with the services sector having the highest share in the economy (\citeauthor{bib31}, \citeyear{bib31}). These high demographic and GDP numbers also correspond with a higher e-commerce activity in Lahore (\citeauthor{bib32}, \citeyear{bib32}).

The data sources for this research mainly include delivery data, road network, and population distribution. The delivery data is gathered from Muller and Phips (M\&P) logistics. M\&P is one the largest logistics operators and provides services all across Pakistan. The road network dataset has been prepared by the National Transport Research Centre (NTRC). NTRC is a government department that deals with transportation research and policy across the country. Worldpop data is used to draw inferences about population distribution. The delivery data set consists of 26,837 `unique Geo-Coordinate' deliveries spread across Lahore. The bulk of these deliveries have been performed in the months of January and February of 2022, over a span of 46 days.
Deliveries with `unique Geo-Coordinate' are defined as those that have distinct pin locations. For instance, a rider may have delivered 10 packages at any particular pin position, which would be represented in the data as 10 different entries. Since each of these deliveries would have the same pin location, they will all be counted as having a single `unique Geo-Coordinate'. The road network acquired by the NTRC includes the following types of roads: primary, secondary, local, motorway, metro, and highway. Motorways and metro are excluded from the network. This is because according to M\&P, most of their deliveries are made using motorbikes which cannot travel on these roads due to local regulations. A direction field is also present in the network for each road segment. Each segment may have one of the following values: North Bound (NB), South Bound (SB), East Bound (EB), West Bound (WB), and None. None roads are two-way direction. 
As a baseline, we clustered the data points according to existing M\&P hubs using the QNEAT3 plugin in QGIS. Points are clustered based on the shortest distance. The delivery point which is closest to a given hub is assigned to that hub. As mentioned previously, QNEAT3 takes road network with directions, hub locations, and delivery locations as input. The Origin-Destination (OD) matrix is computed by the plugin which contains the shortest distance between an origin-destination pair. Distances are determined using the road network. The average distance traveled per delivery from a hub to a delivery location is 7801 meters. It is worth mentioning the delivery model of M\&P at this point. Riders collect their parcels from LH and travel to their respective areas. Therefore, only the first delivery made by riders should be taken into consideration for this study. However, given deliveries are uniformly distributed across Lahore, any delivery can be the first one of the day.

\begin{figure}
\centering
\graphicspath{ {./images/} }
\includegraphics[width = 0.45\columnwidth]{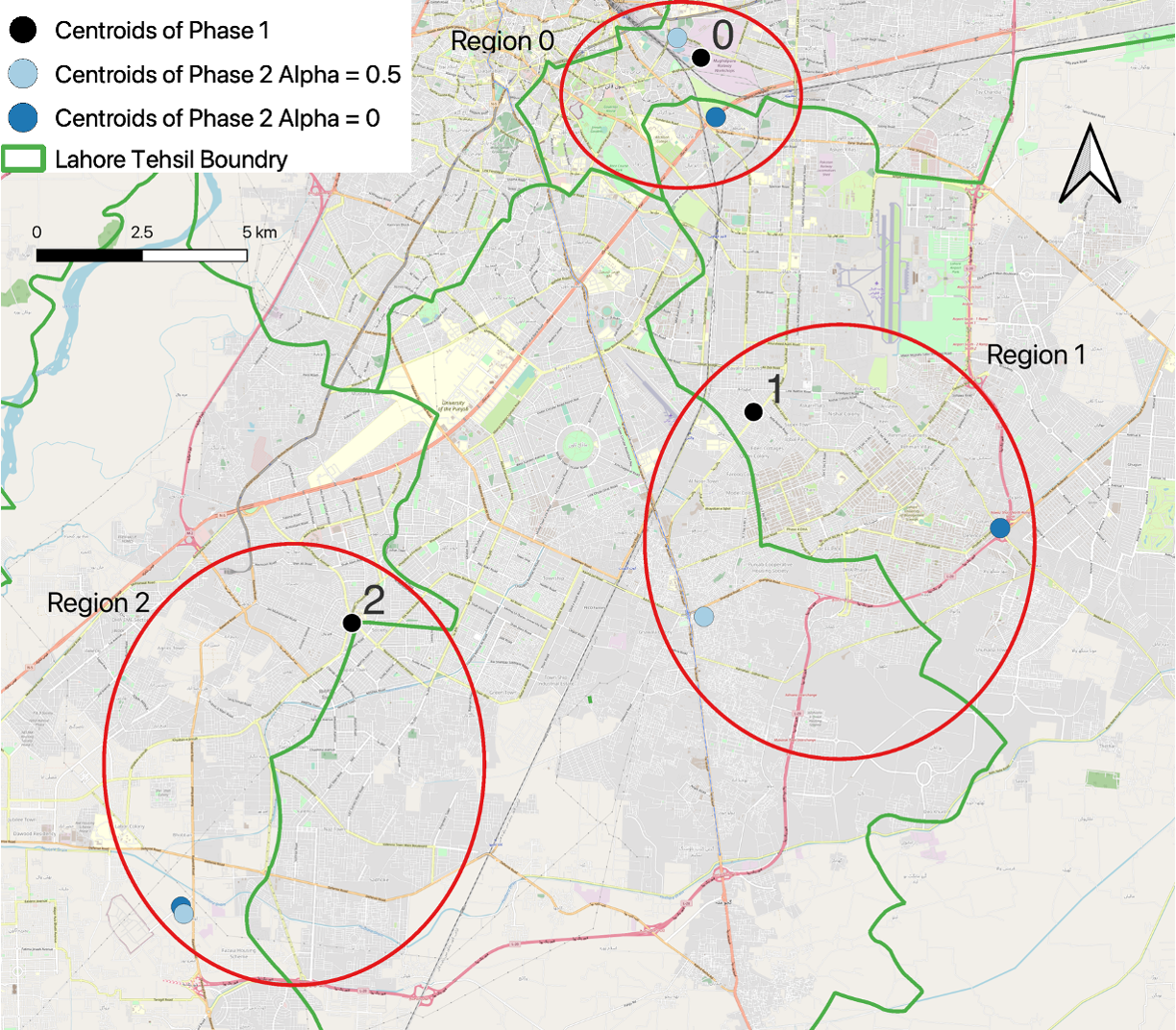}
\caption{Centroids of clusters for phase 1 and phase 2.}
\label{Fig: 2}
\end{figure}
\section{Results}\label{Results}
In phase 1, using M\&P deliveries data, we ran our algorithm for 10 iterations. Current M\&P hubs are used as initial centroids. The cutoff value is 10 meters. The total number of deliveries made at a `unique Geo-Coordinate' are used as weights. For phase 2, the Worldpop dataset is used for clustering where the normalized number of deliveries and normalized ppp are used as weights. Using these weights, the algorithm is run for $\alpha$ values of 0.5 and 0. The cutoff value and number of iterations remain the same. Fig. \ref{Fig: 1} shows the average distance per delivery for each iteration along with convex hulls of clusters formed for phase 1. Fig. \ref{Fig: 1(f)} is a monotonically decreasing graph which means that the algorithm is minimizing average distance during last-mile delivery.
With new clusters and LH in iteration 5, the average distance traveled per delivery decreases from 7801 meters to 6986 meters. The proposed framework reduces the average distance per delivery by 815 meters. The algorithm stops at iteration number 6 since the cutoff criterion is met. Fig. \ref{Fig: 2} shows all the centroids when using only the number of deliveries as a weight (phase-1) and incorporating population as a weight (phase-2). For regions 0 and 1, centroids in phase 2 get closer to the centroids of phase 1 as the value of $\alpha$ increases. For region 2, there is no significant change in the centroids. This shows that for regions 0 and 1, the ``pull" of the number of deliveries as weight is significant as compared to region 2. This effect can also be seen in Table \ref{Table 1}.
\begin{table*} 
\renewcommand{\arraystretch}{2}
\centering
\begin{tabular}{l |l |l}
Centroid phase 1 & Phase 2 $\alpha$ = 0.5 / m & Phase 2 $\alpha$ = 0 / m \\
\hline
0 & 1047 & 2692 \\
1 & 6931 & 8072\\ 
2 & 12324 & 12743\\ 
\end{tabular}
\caption{\label{Table 1} Road network distance to nearest phase 2 centroids from phase 1 centroids.}
\end{table*}
\section{Conclusion}\label{Conclusion}
Previous research focuses on a non-road network-based distance matrix to solve the problem of LH optimization. However, this matrix does not provide a precise measure for distance calculations between points and could lead to misleading results. The most significant effort to address this problem has been to multiply the distance by a constant factor to approximate road network distances. This study intends to bridge the gap in the literature by solving the problem of LH optimization using road network-based distance calculations. A hybrid approach has been developed that utilizes the K-Means algorithm to cluster delivery points. Based on these clusters, a P-Median approach is adopted to find the location of LH. Delivery data points, road network, and population data are used as input. Results show that optimized allocation of LH can significantly reduce the average distance traveled per delivery. The reduction in the average distance would lead to less fuel and time being spent on a particular journey. This would not only positively impact the financial and operational workings of logistics companies, but would also ensure a cleaner environment by cutting greenhouse gas emissions of rider vehicles. This work will inspire industry-standard geospatial solutions which would help alleviate operational concerns of the logistics industry in developing countries and provide a breathing space for the already degrading environment.
A possible area of further research is rider optimization: determining an optimal number of riders to be allocated to a region based on the number of deliveries in that region to reduce operations and management costs. Another area would be to replicate the approach in this study while determining the weights of each possible hub location using various factors. These factors could include CO$_2$ emissions, created job opportunities, noise level, accident risk, and others (\citeauthor{bib35}, \citeyear{bib35}).
\section*{Acknowledgements}
We acknowledge the support provided by M\&P Logistics and NTRC for their active collaboration and for providing us with the necessary data sets used in this research work.


\begin{thebibliography}{}

\bibitem[Ahmed, 2016]{bib32}Ahmed, S.O. (2016). Cracking {E}-commerce 2.0: {W}hitepaper on taking 500,000 merchants online in {P}akistan. Karandaaz.

\bibitem[Agyabeng-Mensah et al., 2020]{bib7}Agyabeng-Mensah, Y., Ahenkorah, E., Afum, E., Dacosta, E., and Tian, Z. (2020). Green warehousing, logistics optimization, social values and ethics and economic performance: the drole of supply chain sustainability. The International Journal of Logistics Management, 31(3), 549-574.

\bibitem[Atayah et al., 2022]{bib2}Atayah, O. F., Dhiaf, M. M., Najaf, K., and Frederico, G. F. (2022). Impact of COVID-19 on financial performance of logistics firms: evidence from G-20 countries. Journal of Global Operations and Strategic Sourcing, 15(2), 172-196.

\bibitem[Black, 2004]{bib22}Black, P. E. (2004). Euclidean {D}istance. Dictionary of Algorithms and Data Structures. \url{https://www.nist.gov/dads/HTML/euclidndstnc.html}

\bibitem[Bozkaya et al., 2010]{bib25}Bozkaya, B., Yanik, S., and Balcisoy, S. (2010). A GIS-based optimization framework for competitive multi-facility location-routing problem. Networks and Spatial Economics, 10, 297-320.

\bibitem[Celik Turkoglu et al., 2020]{bib12}Celik Turkoglu, D., and Erol Genevois, M. (2020). A comparative survey of service facility location problems. Annals of Operations Research, 292, 399-468.

\bibitem[DispatchTrack, 2021]{bib10}DispatchTrack. (2021). Breaking {D}own {L}ast {M}ile {D}elivery {C}osts. \url{https://www.dispatchtrack.com/blog/last-mile-delivery-costs-breakdown}

\bibitem[Dijkstra, 2022]{bib28}Dijkstra, E. W. (2022). A note on two problems in connexion with graphs. Edsger Wybe Dijkstra: His Life, Work, and Legacy (pp. 287-290).

\bibitem[Esnaf et al., 2009]{bib16}Esnaf, Ş., and Küçükdeniz, T. (2009). A fuzzy clustering-based hybrid method for a multi-facility location problem. Journal of Intelligent Manufacturing, 20, 259-265.

\bibitem[Gocer et al., 2022]{bib14}Gocer, F., and Sener, N. (2022). Spherical fuzzy extension of AHP‐ARAS methods integrated with modified k‐means clustering for logistics hub location problem. Expert Systems, 39(2), e12886.

\bibitem[Higgins, 2011]{bib6}Higgins, C. (2011). An exploration of the freight village concept and its applicability to {O}ntario. McMaster Institute For Transportation and Logistics.

\bibitem[ITF, 2019]{bib8}International Transport Forum. (2019). ITF Transport Outlook 2019.

\bibitem[Jin et al., 2017]{bib27}Jin, X., and Han, J. (2017). K-Means clustering. C. Sammut and G. I. Webb (Eds.), Encyclopedia of Machine Learning and Data Mining. (pp. 695-697). Springer.

\bibitem[Kuby et al., 2009]{bib26}Kuby, M., Lines, L., Schultz, R., Xie, Z., Kim, J. G., and Lim, S. (2009). Optimization of hydrogen stations in Florida using the flow-refueling location model. International journal of hydrogen energy, 34(15), 6045-6064.

\bibitem[LCCI, 2017]{bib31}LCCI. (2017). The {S}ize and {G}rowth of the {E}conomy of {L}ahore.

\bibitem[Leblow, 2021]{bib3}Leblow, S. (2021). Worldwide ecommerce continues double-digit growth following pandemic push to online. Emarketer. \url{https://www.emarketer.com/content/worldwide-ecommerce-continues-double-digit-growth-following-pandemic-push-online}

\bibitem[Li et al., 2011]{bib18}Li, Y., Liu, X., and Chen, Y. (2011). Selection of logistics center location using Axiomatic Fuzzy Set and TOPSIS methodology in logistics management. Expert systems with applications, 38(6), 7901-7908.

\bibitem[Mrabti et al., 2022]{bib35}Mrabti, N., Hamani, N., Boulaksil, Y., Gargouri, M. A., and Delahoche, L. (2022). A multi-objective optimization model for the problems of sustainable collaborative hub location and cost sharing. Transportation Research Part E: Logistics and Transportation Review, 164, 102821.

\bibitem[Mulesa et al., 2021]{bib13}Mulesa, O., Mitsa, O., Radivilova, T., Povkhan, I., and Melnyk, O. (2021). Development of a Method to Find the Location of a Logistics Hub. In IT\&I (pp. 263-271).

\bibitem[Önden et al., 2018]{bib17}Önden, İ., Acar, A. Z., and Eldemir, F. (2018). Evaluation of the logistics center locations using a multi-criteria spatial approach. Transport, 33(2), 322-334.

\bibitem[OptimoRoute, 2020]{bib9}OptimoRoute. (2020). What {I}s {L}ast {M}ile {D}elivery? {C}osts \& {H}ow to {O}ptimize. \url{https://optimoroute.com/last-mile-delivery/}

\bibitem[Özceylan et al., 2016]{bib19}Özceylan, E., Erbaş, M., Tolon, M., Kabak, M., and Durğut, T. (2016). Evaluation of freight villages: A GIS-based multi-criteria decision analysis. Computers in Industry, 76, 38-52.

\bibitem[PBS, 2017]{bib29}PBS. (2017). Lahore District. \url{https://www.pbs.gov.pk/census-2017-district-wise/results/053}

\bibitem[Qin et al., 2017]{bib11}Qin, Z., and Gao, Y. (2017). Uncapacitated p-hub location problem with fixed costs and uncertain flows. Journal of Intelligent Manufacturing, 28(3), 705-716.

\bibitem[Shahparvari et al., 2020]{bib5}Shahparvari, S., Nasirian, A., Mohammadi, A., Noori, S., and Chhetri, P. (2020). A GIS-LP integrated approach for the logistics hub location problem. Computers and Industrial Engineering, 146, 106488.

\bibitem[Sharma et al., 2017]{bib15}Sharma, A., and Jalal, A. (2017). Clustering based hybrid approach for facility location problem. Management Science Letters, 7(12), 577-584.

\bibitem[Shezad, 2021]{bib4}Shezad, A. (2021). Pakistan e-commerce platform {D}araz aims to beef up as {A}mazon eyes market. Reuters. \url{https://www.reuters.com/business/retail-consumer/pakistan-e-commerce-platform-daraz-aims-beef-up-amazon-eyes-market-2021-11-25/}

\bibitem[Uluta{\c{s}} et al., 2020]{bib21}Ulutaş, A., Karakuş, C. B., and Topal, A. (2020). Location selection for logistics center with fuzzy SWARA and CoCoSo methods. Journal of Intelligent and Fuzzy Systems, 38(4), 4693-4709.

\bibitem[Wen et al., 2016]{bib23}Wen, R., Yan, W., and Zhang, A. N. (2016, December). Weighted clustering of spatial pattern for optimal logistics hub deployment. In 2016 IEEE International Conference on Big Data (Big Data) (pp. 3792-3797). IEEE.

\bibitem[Worldpop, 2018]{bib34}WorldPop, School of Geography and Environmental Science, University of Southampton; Department of Geography and Geosciences, University of Louisville; Departement de Geographie, Universite de Namur; and Center for International Earth Science Information Network, Columbia University. (2018). Global High Resolution Population Denominators Project - Funded by The Bill and Melinda Gates Foundation (OPP1134076)

\bibitem[Yang et al., 2022]{bib1}Yang, Z., Chen, X., Pan, R., and Yuan, Q. (2022). Exploring location factors of logistics facilities from a spatiotemporal perspective: A case study from Shanghai. Journal of Transport Geography, 100, 103318.

\bibitem[Yang et al., 2020]{bib24}Yang, J., Han, Y., Wang, Y., Jiang, B., Lv, Z., and Song, H. (2020). Optimization of real-time traffic network assignment based on IoT data using DBN and clustering model in smart city. Future Generation Computer Systems, 108, 976-986.

\bibitem[Yang et al., 2016]{bib20}Yang, K., Yang, L., and Gao, Z. (2016). Planning and optimization of intermodal hub-and-spoke network under mixed uncertainty. Transportation Research Part E: Logistics and Transportation Review, 95, 248-266.


































\end{thebibliography}
\end{document}